\documentclass[oneside,english,12pt]{amsart}
\usepackage[latin1]{inputenc}

\makeatletter

\providecommand{\LyX}{L\kern-.1667em\lower.25em\hbox{Y}\kern-.125emX\@}

   \theoremstyle{plain}
   \newtheorem{thm}{Theorem}[section]
   \numberwithin{equation}{section} 
   \numberwithin{figure}{section} 
   \theoremstyle{plain}
   \newtheorem{cor}[thm]{Corollary} 
   \theoremstyle{plain}
   \newtheorem{lem}[thm]{Lemma} 
   \theoremstyle{definition}
   \newtheorem{definition}{Definition}[section]

\usepackage{amsmath}

\newcommand{\image}{\operatorname{im}}
\renewcommand{\Im}{\operatorname{Im}}
\renewcommand{\Re}{\operatorname{Re}}

\makeatother
\usepackage{babel}
\begin{document}

\title{Non-microstates free entropy dimension for groups}

\author{Igor Mineyev$^{\dagger }$ and Dimitri Shlyakhtenko$^{\ddagger
}$}

\begin{abstract}
We show that for any discrete finitely-generated group $G$ and any
self-adjoint $n$-tuple $X_{1},\ldots ,X_{n}$ of generators of the
group algebra $\mathbb{C}G$, Voiculescu's non-microstates free entropy
dimension $\delta ^{*}(X_{1},\ldots ,X_{n})$ is exactly equal to
$\beta _{1}(G)-\beta _{0}(G)+1$, where $\beta _{i}$ are the
$\ell^{2}$-Betti numbers of $G$.
\end{abstract}

\date{December 2003}

\thanks{$\dagger $Research partially supported by NSF CAREER grant
DMS-0228910
and NSF grant DMS-0111298;
the first author also would like to thank the hospitality of IAS,
Princeton, in 2003-04.\\
\indent $\ddagger $Research supported
by the Sloan Foundation and NSF Grant DMS-0102332.}

\maketitle

\section{Introduction.}

In \cite{dvv:entropy2}, using ideas from his theory of free entropy
and free probability, D. Voiculescu has associated to every $n$-tuple
of self-adjoint elements $(X_{1},\ldots ,X_{n})$ in a tracial von
Neumann algebra a number $\delta (X_{1},\ldots ,X_{n})$, which he
called the free entropy dimension of this $n$-tuple. The free entropy
dimension is, very roughly, a kind of asymptotic Minkowski dimension
of the set of $n$-tuples of matrices that approximate the variables
$X_{1},\ldots ,X_{n}$ in non-commutative moments (these are commonly
known as {}``sets of microstates'', see
\cite{dvv:entropysurvey,dvv:entropy3,jung:packing}
for further details).

It is hoped that this number is an invariant of the von Neumann algebra
generated by $X_{1},\ldots ,X_{n}$. While this hope is presently
out of reach in the most interesting cases, this quantity has played
a key role in the solution of several long-standing von Neumann algebra
problems (see e.g. \cite{dvv:entropysurvey} for a survey).

Nonetheless, it is known that a certain technical modification of
$\delta $, $\delta _{0}$ depends only on the algebra generated by
$X_{1},\ldots ,X_{n}$ (and the ambient trace). In particular, if
we start with a discrete finitely-generated group $G$, then $\delta
_{0}$,
evaluated on any set of generators of $G$ gives the same number,
which is an invariant of $G$. This invariant is quite mysterious,
and its exact value is known in only a few cases (such as free products
of abelian groups).

In \cite{dvv:entropy5}, Voiculescu has further introduced a different
approach to free entropy and free entropy dimension, based on the
theory of free Hilbert transform. This {}``microstates-free'' approach
has resulted in two definitions of {}``non-microstates'' free entropy
dimension-like quantities, $\delta ^{*}$ and $\delta ^{\star }$.
While it is suspected that $\delta ^{*}=\delta ^{\star }$, we only
know that always $\delta ^{\star }\geq \delta ^{*}$. By a deep result
of Biane, Capitaine and Guionnet
\cite{guionnet-biane-capitaine:largedeviations},
$\delta ^{*}\geq \delta $.

Much less is known about $\delta ^{*}$ than about $\delta $; in
all the known cases they assume the same value, although this statement
speaks more for the small number of cases in which the value of both
is known than for the existence of a general strategy to prove that
they are the same for some class of $n$-tuples. Only recently have
there been any non-trivial computations of $\delta ^{*}$
(\cite{aagaard:freedim},\cite{shlyakht:qdim}).

Let $G$ be a finitely generated
discrete group, and let $\mathbb{C}G$ be its
group algebra, endowed with the involution $
(\sum_\gamma \alpha_\gamma \gamma )^* = \sum_\gamma
\bar\alpha_\gamma\gamma^{-1}$
and the tracial linear functional $\tau(\sum_\gamma \alpha_\gamma
\gamma)
=\alpha_e$.  Let $X_1,\ldots,X_n$ be any generators of this algebra,
which are self-adjoint (e.g., if $\gamma_1,\ldots,\gamma_m$ are
generators of~$G$
one could take $n=2m$ and
$X_j = \gamma_j + \gamma_j^{-1}$, $1\leq j\leq m$,
$X_j = -i( \gamma_{j-m} - \gamma_{j-m}^{-1}), j=m+1,\ldots,2m$.

Recently, in \cite{connes-shlyakht:l2betti} A. Connes and the second
author have proved that
\begin{equation}
\delta ^{*}(X_{1},\ldots ,X_{n})\leq \delta ^{\star }(X_{1},\ldots
,X_{n})\leq \beta _{1}(G)-\beta
_{0}(G)+1,\label{eq:deltastarvsbettis}\end{equation}
   where $\beta _{i}(G)$ are Atiyah's $\ell^{2}$-Betti numbers of the
group $G$ (see \cite{atiyah-L2,cheeger-gromov:l2,luck:book}). The
appearance of $\ell^{2}$-invariants of $G$ in connection with free
entropy dimension has been conjectured by specialists
ever since the
fundamental work of Gaboriau \cite{gaboriau:ell2,gaboriau:cost}.
Nonetheless, this connection remains quite surprising
to us, since
free entropy dimension relies on the notion of free Brownian motion,
while $\ell^2$-Betti numbers are homological in nature,
and it is hard to say why the two must have anything in common.

The main result of this paper is that in fact equality holds: we prove
that\[
\delta ^{*}(X_{1},\ldots ,X_{n})=\delta ^{\star }(X_{1},\ldots
,X_{n})=\beta _{1}(G)-\beta _{0}(G)+1,\]
for any finitely-generated group $G$ and any set of self-adjoints
$X_{1},\ldots ,X_{n}\in\mathbb{C}G$ generating $\mathbb{C}G$. In
particular, we
conclude that in this case, $\delta ^{*}=\delta ^{\star }$, and both
are algebraic invariants.

The main technical tool is a result showing that arbitrary $\ell ^{2}$
1-coboundaries on the Cayley graph of $G$ can be approximated in $\ell
^{2}$
norm by coboundaries of the form $\delta g$, where
$g\in \ell ^{\infty }(G)$. This result holds more
generally for arbitrary graphs, and for $\ell ^{2}$ replaced by $\ell
^{p}$,
$1\leq p<\infty $.

Using this result, we utilize a lower estimate for non-microstates
free entropy dimension from \cite{shlyakht:qdim}, which combined
with (\ref{eq:deltastarvsbettis}) gives the main result.

\section*{Notations.}
Throughout this paper, $G$ will denote a finitely generated discrete
group.  We write $\ell^2(G)$ for the Hilbert space of square-summable
functions on $G$.  We denote by $\lambda$ and $\rho$ the left and
right regular representation of $G$ on $\ell^2(G)$, and by $L(G)$ the
group
von Neumann algebra, which is the weak operator topology closure
of the linear span of the image $\lambda(G)$ viewed as subalgebra
of the algebra of bounded operators $B(\ell^2(G))$.  By $\tau$ we shall
always denote the von Neumann trace on $L(G)$ given by $\tau(x)
=\langle x \delta_e,\delta_e\rangle$, where $\delta_e$ is the delta
function at the identity of $G$.   The restriction of $\tau$ to
$\lambda(\mathbb{C}G)$ is the canonical group trace
on the group algebra determined by linearity
and the condition $\tau(g)=1$ if $g=e$ and $\tau(g)=0$ if $g\neq e$.

The letter $M$ will denote a general von Neumann algebra with a normal
(i.e. weak-operator continuous)
tracial state $\tau: M\to \mathbb C$, $\tau(xy)=\tau(yx)$.  The von
Neumann
algebra $M$ acts by left and right multiplication
on the Hilbert space $L^2(M)$, which is the completion
of $M$ in the norm $\Vert m\Vert_2 = \tau(m^*m)^{1/2}$.
    In the case that $M=L(G)$, $L^2(M)=\ell^2(G)$, and the
left and right actions of $L(G)$ on this space extend the left and right
actions of $G$.  We will denote by $M^o$ the opposite von Neumann
algebra.
The letter $J$ will denote the anti-linear
Tomita conjugation operator $J:L^2(M)\to
L^2(M)$ extending $J(m)= m^*$.  The operator $J$ satisfies the property
that for $x\in M$ and $\xi\in L^2(M)$, $JxJ \ \xi = \xi \ x^*$, i.e., it
switches the right and left actions of $M$.  In particular, for any
$x\in M$, $JxJ$ commutes with $M$.

If $H\subset L^2(M)^{\oplus n}$ is a closed $M$-submodule of a multiple
of the left module $L^2(M)$, we denote by $\dim_{M} H$ its Murray-von
Neumann
dimension.  This dimension satisfies the usual monotonicity and
additivity
properties (see Chapter X in \cite{MvN} (esp. Theorem X on p. 182), or, 
for a more accessible introduction, \cite{GHJ,cheeger-gromov:l2}).

We denote by $B(L^2(M))$ the space of all bounded linear operators on
$L^2(M)$.
Finally, we will denote by $HS$ the space of Hilbert-Schmidt operators
$T:L^2(M)\to L^2(M)$, i.e. the operators $T\in B(L^2(M))$ for which the
norm
$\Vert T\Vert_{HS} = \operatorname{Tr}(T^*T)$ is finite.  $HS$ is a
Hilbert space with the inner product $\langle
T,S\rangle=\operatorname{Tr}
(TS^*)$.  $HS$ can be identified with the Hilbert space tensor product
$L^2(M)\bar\otimes L^2(M^o)$ by the map $m\otimes n^o\mapsto m P_1 n$,
where $P_1\in HS$ denotes the rank one projection onto the vector
$1\in M\subset L^2(M)$, $m\in M$ and $n^o\in M^o$.
By definition, the von Neumann algebra tensor product $M\bar\otimes M^o$
acts on the tensor product Hilbert space $L^2(M)\bar\otimes L^2(M^o)$
and
thus on $HS$.
The action of the algebraic tensor product
$M\otimes M^o\subset M\bar\otimes M^o$ on $HS$ is explicitly given by
   $(m\otimes n^o)\cdot T= mTn$ (composition of operators
on $L^2(M)$), for $T\in HS$, $m\in M$
and $n^o\in M^o$.

\section{Approximation of $\ell ^{p}$-summable 1-coboundaries on
graphs.}

Let $\mathcal{G}$ be a graph. $C^{i}(\mathcal{G},\mathbb{R})$ will
denote the set of real $i$-cochains on~$\mathcal{G}$, without any
assumptions on their support.  For each $1\leq p\leq \infty $, let
$C_{(p)}^{i}(\mathcal{G},\mathbb{R})$ be the set of elements in
$C^{i}(\mathcal{G},\mathbb{R})$ which
have finite $\ell ^{p}$ norm. Finally, let
$\delta :C^{0}(\mathcal{G},\mathbb{R})\to C^{1}(\mathcal{G},\mathbb{R})$
be the coboundary map.

\begin{thm}
\label{approximation} Let $\mathcal{G}$ be an arbitrary graph,
$p\in [1,\infty )$ and $f\in C^{0}(\mathcal{G},\mathbb{R})$
be such that $\delta f\in C_{(p)}^{1}(\mathcal{G},\mathbb{R})$. Then
for each $\varepsilon >0$ there exists $g\in C_{(\infty
)}^{0}(\mathcal{G},\mathbb{R})$
such that $\Vert \delta f-\delta g\Vert _{p}<\varepsilon $. In
particular,
$\delta g\in C_{(p)}^{1}(\mathcal{G},\mathbb{R})$.

The same result holds in the complex-valued case.
\end{thm}
\begin{proof}
Let $\Sigma _{i}$ denote the set of $i$-simplices in~$\mathcal{G}$.
We are given a function $f:\Sigma _{0}\to \mathbb{R}$ such that $\delta
f:\Sigma _{1}\to \mathbb{R}$
is $\ell ^{p}$-summable.

Fix some $t\in [0,\infty )$, denote $U_{t}=f^{-1}([-t,t])\subseteq
\Sigma _{0}$
and for $x\in \Sigma _{0}$, \[
f_{t}(x)=\begin{cases} -t & \textrm{if}\  f(x)\in (-\infty ,-t)\\
   f(x) & \textrm{if}\  f(x)\in [-t,t]\\
   t & \textrm{if}\  f(x)\in (t,\infty ).\\ \end{cases}
\]
   Obviously, $|f_{t}(x)|\le \min \{|f(x)|,t\}\le t$, so in particular
$f_{t}\in C_{(\infty )}^{0}(\mathcal{G},\mathbb{R})$ for each $t$.

Let $\delta U_{t}$ be the set of all edges in $\mathcal{G}$ all
of whose incident vertices are in $U_{t}$. We have \[
\bigcup _{t\in [0,\infty )}U_{t}=\Sigma _{0},\]
   and therefore \begin{equation}
\bigcup _{t\in [0,\infty )}\delta U_{t}=\Sigma
_{1},\label{cover}\end{equation}
   where $\{\delta U_{t}\}$ is an increasing sequence of sets.

Since $f$ and $f_{t}$ coincide on $U_{t}$, then $\delta f$ and
$\delta f_{t}$ coincide on $\delta U_{t}$, that is \begin{equation}
\operatorname {supp}(\delta f-\delta f_{t})\subseteq \Sigma
_{1}\setminus \delta U_{t}.\label{support}\end{equation}
We need the following lemma.

\begin{lem}
\label{boundary-bound} With the above notations, \[
|\delta f_{t}(e)|\le |\delta f(e)|\textrm{ for all }t\in [0,\infty
)\textrm{ and }e\in \Sigma _{1}.\]

\end{lem}
\begin{proof}
Since $\delta f$ and $\delta f_{t}$ coincide on $\delta U_{t}$,
it only remains to show the inequality when $e\in \Sigma _{1}\setminus
\delta U_{t}$,
that is when the edge $e$ is incident to a vertex $x$ in~$\Sigma
_{0}\setminus U_{t}$.
By the definition of $U_{t}$ this means that $f(x)\in (-\infty ,-t)\cup
(t,\infty )$.
We can assume $f(x)\in (t,\infty )$, the opposite case can be done
similarly. Let $x'$ be the other incident vertex of~$e$. There
are three obvious cases to consider for $x'$, and we use the definition
of $f_{t}$ in each case.

If $f(x')\in (t,\infty )$ then \[
|\delta f_{t}(e)|=|f_{t}(x)-f_{t}(x')|=|t-t|=0\le |\delta f(e)|.\]
   If $f(x')\in [-t,t]$ then \begin{eqnarray*}
   &  & |\delta f_{t}(e)|=|f_{t}(x)-f_{t}(x')|=|t-f(x')|=t-f(x')\\
   &  & \le f(x)-f(x')=|f(x)-f(x')|=|\delta f(e)|.
\end{eqnarray*}
   If $f(x')\in (-\infty ,-t)$ then \begin{eqnarray*}
   &  & |\delta f_{t}(e)|=|f_{t}(x)-f_{t}(x')|=|t-(-t)|=t+t\\
   &  & \le f(x)-f(x')=|f(x)-f(x')|=|\delta f(e)|.
\end{eqnarray*}
   This finishes the proof of the lemma.
\end{proof}
Now we can finish the proof of Theorem~\ref{approximation}.
Since $\delta f$ is $\ell ^{p}$-summable, given any $\varepsilon >0$,
(\ref{cover}) guarantees the existence of $t\in [0,\infty )$ such
that \[
\Vert \delta f|_{\Sigma _{1}\setminus \delta U_{t}}\Vert
_{p}<\varepsilon /2,\]
   then by Lemma~\ref{boundary-bound}, \[
\Vert \delta f_{t}|_{\Sigma _{1}\setminus \delta U_{t}}\Vert _{p}\le
\Vert \delta f|_{\Sigma _{1}\setminus \delta U_{t}}\Vert
_{p}<\varepsilon /2,\]
   so (\ref{support}) implies that \[
\Vert \delta f-\delta f_{t}\Vert _{p}=\Vert (\delta f-\delta
f_{t})|_{\Sigma _{1}\setminus \delta U_{t}}\Vert _{p}=\Vert \delta
f|_{\Sigma _{1}\setminus \delta U_{t}}-\delta f_{t}|_{\Sigma
_{1}\setminus \delta U_{t}}\Vert _{p}\le \varepsilon /2+\varepsilon
/2=\varepsilon .\]

Setting $g=f_{t}$ completes the proof of Theorem~\ref{approximation}
in the real case. The complex case is obtained by separately
approximating
the real and imaginary parts of $\delta f$.
\end{proof}

\section{$\ell^{2}$-Betti Numbers.}

\subsection{$\ell ^{2}$-Betti Numbers for Groups.}
The notion of $\ell^2$-Betti numbers for groups goes back to Atiyah
\cite{atiyah-L2} and Cheeger and Gromov \cite{cheeger-gromov:l2}.
We refer the reader to the book \cite{luck:book} for more details and
only sketch the construction here.

Assume that
the group $G$ acts freely on a CW-complex $X$, and that the complex $X$
is
``co-finite'' (i.e., for each dimension $i$ there is a finite number of
$i$-cells
in $X$, so that every other $i$-cell in $X$ can be obtained from one of
them
by the group action).  Let $C_i^{(2)}(X,\mathbb {C})$
denote the complex Hilbert space whose
orthonormal basis is formed by the $i$-cells of the complex $X$.  Then
$G$ acts on $C_i^{(2)}(X)$; this action of course extends to a
representation
of the group algebra $\mathbb{C}G$ of $G$ on this Hilbert space.
This representation is contained in a multiple
of the
left regular representation,
and hence the action of $\mathbb{C}G$
extends by continuity to an action of the group von Neumann
algebra $L(G)$.

Thus one can speak of the Murray-von Neumann
dimension of any closed $G$-invariant subspace of the Hilbert space
$C_i^{(2)}(X)$.

The boundary maps $\partial_i$ of the complex $X$ extend to continuous
linear operators $\hat\partial_i: C_i^{(2)}(X,\mathbb{C}) \to
C_{i-1}^{(2)}(X,
\mathbb{C})$.

The reduced $\ell^2$-homology of the complex $X$ is defined to be
the sequence of Hilbert spaces $$
H^{(2)}_k(X)=\ker{\hat\partial_k}/ \overline{\image \partial_{k+1}},$$
where closure
is taken with respect to the Hilbert space norm
(the closure of $\image\partial_{k+1}$ is the same as that of
$\image\hat\partial_{k+1}$).
Note that $H^{(2)}_k(X)$ can be
thought of as the orthogonal complement of $\image\partial_{k+1}$ inside
$\ker \hat\partial_k \subset C_k^{(2)}(X,\mathbb{C})$.
   Thus one can consider its
Murray-von Neumann dimension, which is exactly the $k$-th $\ell^2$-Betti
number of $(X,G)$:
$$
\beta_k(X,G)=\dim_{L(G)} H_k^{(2)}(X).
$$
In the case that $X$ is not co-finite, one writes $X$ as an increasing
union of
co-finite $G$-invariant subcomplexes $X_n$, $n=1,2,\ldots$.  In that
case
the $\ell^2$-Betti numbers can be computed as the following limits
\begin{equation}
\label{Betti-k}
\beta_k(X,G) = \sup_{n} \inf_{m\geq n} \dim_{L(G)} \frac{
    \ker \hat\partial_k : C_k^{(2)}(X_n,\mathbb{C})\to
C_{k-1}^{(2)}(X_n,\mathbb{C})}
    {\overline{\left(
    \image \hat\partial_{k+1}:C_{k+1}^{(2)}(X_m,\mathbb{C})\to
C_k^{(2)}(X_m,\mathbb{C})
    \right) }\cap C_k^{(2)}(X_n,\mathbb{C})}.
\end{equation}
(closure in Hilbert space norm, see \cite{cheeger-gromov:l2}).

The main point of interest for us is the fact that if the CW-complex $X$
is $n$-connected, then the first $n+1$ $\ell^2$-Betti numbers
$\beta_0(X,G),
\beta_1(X,G),\ldots,\beta_{n}(X,G)$
are independent of $X$ and are invariants of the group $G$.
In this case,
they are referred to as the $\ell^2$-Betti numbers of the group $G$.

\subsection{Zeroth and First $\ell^2$-Betti numbers for
finitely-generated
groups.}

If (as we are in the present paper) one is only interested in the
zeroth and first $\ell^2$-Betti numbers of a finitely generated group
$G$, then one can make an explicit choice of a one-connected CW-complex
that can be used to compute the first two $\ell^2$-Betti numbers.

Let $\mathcal{G}$ denote the Cayley graph of $G$ with respect to the
set of generators $g_{1},\ldots ,g_{n}$. Then $G$ acts on
$\mathcal{G}$ by left translation. We view $\mathcal{G}$ as a
CW-complex, whose $1$-cells are the edges of $\mathcal{G}$ and whose
0-cells are the vertices of $\mathcal{G}$. There exists a
simply-connected CW-complex $X$, whose $1$-skeleton is $\mathcal{G}$;
it is obtained from $\mathcal{G}$ by gluing in a single 2-cell for
each non-trivial loop in $\mathcal{G}$.

The action of $G$ on the CW-complex $X$ need not be co-finite
(although it is co-finite when restricted to the $1$-skeleton, since
the group $G$ is finitely-generated).  However, one can write $X$ as a
union of $X_m$, $m=1,2,\ldots$, where $X_m$ are $G$-invariant
subcomplexes of $X$, having $\mathcal{G}$ as their $1$-skeletons, and
with the property that each $X_m$ is co-finite.  Indeed, one could
just enumerate all of the $2$-cells used in the construction of $X$,
and for each $m$, let $X_m$ be the space arising after the first $m$
$2$-cells, together with all of their $G$-translates, are glued to
$\mathcal{G}$.

We consider the spaces of $i$-cells of $X_m$ as subsets 
$C_{i}(X_m)\subset C_{i}(X)$.
Let us denote by $C_{i}^{(2)}(X_m,\mathbb{C})$ the completion of the 
space
$C_{i}(X_m,\mathbb{C})$ with respect to $\ell ^{2}$-norm.
Let $\partial _{i}:C_{i}(X_m,\mathbb{C})\to C_{i-1}(X_m,\mathbb{C})$ be
the boundary map and
${\hat{\partial}}_{1}:C_{1}^{(2)}(X_m,\mathbb{C})\to
C_{0}^{(2)}(X_m,\mathbb{C})$, $i=1,2$, be its continuous extension.

In this case \cite{cheeger-gromov:l2,bekka-valette:l2cohomology} the
first two $\ell ^{2}$-Betti numbers of $G$ are defined as the following
Murray-von Neumann dimensions over the group von Neumann algebra $L(G)$
of $G$:\begin{equation*}
\beta_{1}(G)  =  \dim_{L(G)}H_{1}^{(2)}(X),\qquad
\beta_{0}(G)  =  \dim_{L(G)}H_{0}^{(2)}(X).
\end{equation*}

Then we have by additivity of dimension
and by~(\ref{Betti-k}),

\begin{eqnarray}\label{beta0}
\nonumber&&\beta_{1}(G)
    =\inf_{m\ge 1} \dim_{L(G)} \frac{
    \ker \hat\partial_1 : C_1^{(2)}(X_1,\mathbb{C})\to
C_{0}^{(2)}(X_1,\mathbb{C})}
    {\overline{
    \image \hat\partial_2:C_{2}^{(2)}(X_m,\mathbb{C})\to
C_1^{(2)}(X_m,\mathbb{C})}}
    =\\
\nonumber&&=\inf_{m}\left( \dim_{L(G)}\ker\hat{\partial}_{1}
-\dim_{L(G)} \overline{\hat\partial_{2}(C_2(X_m))}\right),
   \\
&&\beta_{0}(G)  =1-\dim _{L(G)}\overline{\image \hat{\partial}_{1}}.
\end{eqnarray}

Note that $\overline{\partial_{2}(C_2(X_m))}$, $m=1,2,\ldots$ are
increasing $L(G)$-submodules of a finite-dimensional $L(G)$-module
$\ker \hat{\partial}_{1}$.  Thus
$$\inf_{m} \left(\dim_{L(G)}\ker\hat{\partial}_{1} -\dim_{L(G)}
\overline{\hat\partial_{2}(C_2(X_m))}\right) = \dim _{L(G)}\ker
\hat{\partial}_{1}-\dim _{L(G)}\overline{ \hat\partial _{2}(C_2(X))}.$$
Since $X$ is
simply-connected,  $\image \partial_{2}= \ker \partial_{1}$ and their
$\ell^2$-closures inside
$C_1^{(2)}$ coincide with the closure of the space
$\image\hat\partial_{2}$.
   Thus
\begin{equation}\label{beta1}
   \beta_1(G)=\dim _{L(G)}\ker
\hat{\partial}_{1}-\dim _{L(G)}\overline{\ker \partial _{1}}.
\end{equation}

Denote by $C^{i}(X,\mathbb{C})$ the space
of all cochains on $X$, i.e., the algebraic dual of
$C_{i}(X,\mathbb{C})$, and by
$$\delta: C^{i}(X,\mathbb{C})\to C^{i+1}(X,\mathbb{C})$$ the
coboundary map.
Let $C_{(2)}^{i}(X,\mathbb{C})$ be the space of all $\ell
^{2}$-summable $i$-cochains on $X$.
Then by duality,
\begin{equation}
\label{kerC}
\overline{\ker \partial _{1}}
    =  \overline{\image\partial_{2}} = \{c\in
C_{(2)}^{1}(X,\mathbb{C}):\delta c=0\}^{\perp }\subset
C_{1}^{(2)}(X,\mathbb{C}).
\end{equation}
Here we identify both
$C_{(2)}^{1}(X,\mathbb{C})$ and $C^{(2)}_{1}(X,\mathbb{C})$ with
$\ell^2(\Sigma_1)$, $\Sigma_1$ being the set of 1-simplices in~$X$,
and all the closures and orthogonal complements are taken in
$\ell^2(\Sigma_1)$.

The first cohomology of the complex
$C^{*}(X,\mathbb{C})$ vanishes, since $X$ is simply-connected.

Therefore if $c\in
C_{(2)}^{1}(X,\mathbb{C})$ satisfies $\delta c=0$, then $c=\delta f$
for some $f\in C^{0}(X,\mathbb{C})$. Thus by~(\ref{kerC}),
\[ \overline{\ker
    \partial_{1}}=\Big(\delta
(C^{0}(X,\mathbb{C}))\cap C_{(2)}^{1}(X,\mathbb{C})\Big)^{\perp }\subset
C_{1}^{(2)}(X,\mathbb{C}).\] Theorem \ref{approximation} says that
$$ \delta(C^{0}(X,\mathbb{C}))\cap C_{(2)}^{1}(X,\mathbb{C})\subseteq
\overline{\delta (C_{(\infty )}^{0}(X,\mathbb{C}))\cap
C_{(2)}^{1}(X,\mathbb{C})},$$
so we get
the following corollary:

\begin{cor}
\label{cor:approxCohomology} The closure of $\image\partial _{2}$
equals\[
\overline{\ker \partial _{1}}=
\Big(\delta (C_{(\infty )}^{0}(X,\mathbb{C}))\cap
C_{(2)}^{1}(X,\mathbb{C})\Big)^{\perp }\subset
C_{1}^{(2)}(X,\mathbb{C}).\]

\end{cor}

\begin{lem}
\label{lem:delta2} Let $\delta ^{(2)}(G)=\beta _{1}(G)-\beta
_{0}(G)+1$.
Then\begin{equation*}
\delta ^{(2)}(G) =  n-\dim _{L(G)}\overline{\ker \partial _{1}}
   =  \dim _{L(G)}\overline{\left(\delta (C_{(\infty
)}^{0}(X,\mathbb{C}))\cap C_{(2)}^{1}(X,\mathbb{C})\right)}.
\end{equation*}

\end{lem}
\begin{proof}
We have by (\ref{beta1}) and (\ref{beta0})
\begin{eqnarray*}
\beta _{1}(G)-\beta _{0}(G)+1 & = & \dim _{L(G)}\ker \hat{\partial}
_{1}-\dim _{L(G)}\overline{\ker \partial _{1}}-1+\dim
_{L(G)}\overline{\image \hat{\partial}_{1}}+1\\
   & = & \dim _{L(G)}\ker \hat{\partial}_{1}+\dim _{L(G)}\overline{\image
     \hat{\partial}_{1}}-\dim _{L(G)}\overline{\ker \partial _{1}}\\
   & = & \dim _{L(G)}C_{1}^{(2)}(X;\mathbb{C})-\dim
   _{L(G)}\overline{\ker \partial _{1}},
\end{eqnarray*}
the last equality by additivity of Murray-von Neumann dimension.
But $C_{1}^{(2)}(X;\mathbb{C})\cong (\ell ^{2}(G))^{\oplus n}$, so
that\[
\delta ^{(2)}(G)=n-\dim _{L(G)}\overline{\ker \partial _{1}
    }=\dim _{L(G)}(\ker \partial _{1})^{\perp }.\]
It remains to apply Corollary \ref{cor:approxCohomology}.
\end{proof}

\subsection{$\Delta$ and $L^{2}$-homology of algebras.}

Let $(M,\tau )$ be a tracial von Neumann algebra,
and let $X_{1},\ldots ,X_{n}\in M$
be a self-adjoint set of elements (i.e., we assume that for each $i$,
there is a $j$ so that $X_{i}^{*}=X_{j}$). Let $HS$ be the space
of Hilbert-Schmidt operators on the Hilbert space $L^{2}(M,\tau )$.

Let $J:L^{2}(M,\tau )\to L^{2}(M,\tau )$ be the anti-linear Tomita
conjugation
operator (see notations).
Then $JMJ$ is exactly the commutant of $M$ in $B(L^2(M))$.

We view $HS$ as a bimodule over $M$ using the action\[
(m_{1}\otimes m_{2}^{o})\cdot T=m_{1}Tm_{2},\qquad m_{1},m_{2}\in
M,\quad T\in H
S.\]
Note that since $HS\cong L^2(M,\tau)\bar\otimes L^2(M,\tau)^o\cong
L^2(M\bar\otimes M^o)$,
the action of
$M\otimes M^o$ on $HS$ extends by continuity to the action of the
von Neumann algebra $M\bar\otimes M^o$, which
is exactly the left multiplication action of $M\bar\otimes M^o$ on
$L^2(M\bar\otimes M^o)$.  In particular, if $H$ is any
$M,M$-sub-bimodule
of $HS$, which is closed in the Hilbert-Schmidt norm, then
it is a module over $M\bar\otimes M^o$; in particular, the
Murray-von Neumann dimension of $H$ over $M\bar\otimes M^o$ makes sense.

Some of the main ideas of the approach
to $L^2$ homology of algebras
in \cite{connes-shlyakht:l2betti}, when particularized to the case
of the first Betti number, can be summarized
in the following (well-known) table, giving a dictionary between
group and von Neumann algebra terms (here $[X,Y] = XY-YX$ denotes the 
commutator
of $X$ and $Y$):

\bigskip
\begin{center}
\begin{tabular}{|l|l|}
\hline
\hfil Group $G$ \hfill &
\hfil von Neumann algebra $M$ \hfill \\ \hline
$\ell^2(G)$ as a group module & $HS$ as an $M,M$-bimodule  \\[7pt]
$g_1,\ldots,g_n$ generators of $G$ &
\begin{minipage}[t]{200pt} $X_j=\lambda_{g_j}$, $j=1,\ldots,n$
in the left regular representation  $\lambda$ of $G$\\
   \end{minipage} \\[20pt]
$\ell^\infty(G)$ & $B(L^2(M,\tau))$ \\[7pt]
Function $f$ on $G$
& Operator $m_f$ of multiplication by $f$ \\[7pt]
\begin{minipage}[t]{200pt}
$\delta f = (\rho_{g_1}(f)-f,\ldots,\rho_{g_n}(f)-f)
\in C^1_{(2)}(\mathcal{G})\cong \ell^2(G)^{\oplus n}$
with $f\in \ell^\infty(G)$ ($\rho$ is the right regular representation)
\end{minipage}&
\begin{minipage}[t]{200pt}
$([D,JX_1J],\ldots,[D,JX_nJ])\in HS^n$ for $D\in B(L^2(M,\tau))$.
(see equation (\ref{eq:groupvsM}) and also Lemma \ref{lem:groupvsM}).
\end{minipage} \\ \hline
\end{tabular}
\end{center}
\bigskip
Here $[\cdot,\cdot]$ denotes the commutator in $B(L^2(M))$.

Following the ideas presented in the table above and
\cite[Corollary 2.12]{shlyakht:qdim}
(we caution the reader that the roles of $M$ and $JMJ$ are switched
in the present paper compared to \cite{shlyakht:qdim}), consider the
set\[
H_{0}(X_{1},\ldots ,X_{n})=\{(\Xi _{1},\ldots ,\Xi _{n})\in
HS^{n}:\exists D\in B(L^{2}(M))\textrm{ s}.\textrm{t. }\Xi _{j}=[D,J
X_{j} J ]\ \forall j\}.\]
Then $H_{0}$ is an $M,M$-bimodule.

   \begin{definition}
   Let\[
\underline{\Delta }(X_{1},\ldots ,X_{n})=\dim _{M\bar{\otimes
}M^{o}}\overline{H_{0}(X_{1},\ldots ,X_{n})},\]
where the closure is taken in the Hilbert-Schmidt topology on $HS$.
   \end{definition}

The quantity $\underline{\Delta}$ has appeared in \cite{shlyakht:qdim}
in
connection with some technical estimates on free entropy dimension.
As we shall see later in Lemma \ref{lemma:delbtaBarvsdelta2}
(and as is apparent from our table of analogies),
the space $H_0(X_1,\ldots,X_n)$ is the von Neumann algebra analog of
the space $$
\{ c \in C_{(2)}^1 (\mathcal{G}) : c = \delta f \textrm{ for some }
   f \in \ell^\infty(G) \}.$$

The proof of the following Lemma was inspired by the work of Bekka
and Valette \cite{bekka-valette:l2cohomology}.

\begin{lem}
Assume that $X_{1},\ldots ,X_{n}$ generate $M$ as a von Neumann
algebra. Then $\underline{\Delta }(X_{1},\ldots ,X_{n})$ depends
only on the algebra $\mathbb{C}(X_{1},\ldots ,X_{n})$ generated by
$X_{1},\ldots ,X_{n}$ and the trace $\tau $.
\end{lem}
\begin{proof}
For $D\in B(L^{2}(M))$, define a Hilbert space seminorm by
\begin{equation}\label{seminorm-def}
\Vert D\Vert _{X_{1},\ldots ,X_{n}}=\left(\sum_{j=1}^n \Vert [D,JX_{j}
J]\Vert _{HS}^{2}\right)^{1/2}.
\end{equation}
Let $\tilde{D}(X_{1},\ldots ,X_{n})=\{D:\Vert D\Vert _{X_{1},\ldots
,X_{n}}<\infty \}$,
and let $D_{0}(X_{1},\ldots ,X_{n})$ be the Hilbert space obtained
from $\tilde{D}(X_{1},\ldots ,X_{n})$ after separation
   and completion.
Endow $D_{0}(X_{1},\ldots ,X_{n})$ with the $M,M$-bimodule structure
coming from the action $(m\otimes n^{o})\cdot D=mDn$. Then the
map\[
D\mapsto ([D,JX_{1}J],\ldots ,[D,JX_{n}J])\]
descends and extends to an $M\bar{\otimes }M^{o}$-module isomorphism
of $D_{0}(X_{1},\ldots ,X_{n})$ with the Hilbert-Schmidt completion
of $H_{0}(X_{1},\ldots ,X_{n})$.

Let $Y_{1},\ldots ,Y_{m}\in \mathbb{C}(X_{1},\ldots ,X_{n})$.
By the definition of the seminorm in~(\ref{seminorm-def})
we clearly have\[
\Vert D\Vert _{X_{1},\ldots ,X_{n}}\leq \Vert D\Vert _{X_{1},\ldots
,X_{n},Y_{1},\ldots ,Y_{m}}.\]
Also, since each $Y_{j}$ is a polynomial in $X_{1},\ldots ,X_{n}$,
   $\Vert [D,JY_{j}J]\Vert _{HS}\leq C_{j}\Vert D\Vert _{X_{1},\ldots
,X_{n}}$
for some constants $C_{1},\ldots ,C_{m}$. It follows that the norms
$\Vert \cdot \Vert _{X_{1},\ldots ,X_{n}}$ and $\Vert \cdot \Vert
_{X_{1},\ldots ,X_{n},Y_{1},\ldots ,Y_{m}}$
are equivalent. Thus the Hilbert space completions of
$H_{0}(X_{1},\ldots ,X_{n})$
and $H_{0}(X_{1},\ldots ,X_{n},Y_{1},\ldots ,Y_{m})$
are isomorphic as $M\bar{\otimes }M^{o}$-modules. Thus
\begin{equation}\label{eqn:canAddGenerators}
\underline{\Delta }(X_{1},\ldots ,X_{n})=\underline{\Delta
}(X_{1},\ldots ,X_{n},Y_{1},\ldots ,Y_{m}).
\end{equation}
If $Y_{1},\ldots ,Y_{m}$ generate $\mathbb{C}(X_{1},\ldots ,X_{n})$,
then  by (\ref{eqn:canAddGenerators}) \[
\underline{\Delta }(Y_{1},\ldots ,Y_{n})=\underline{\Delta
}(X_{1},\ldots ,X_{n},Y_{1},\ldots ,Y_{m})=\underline{\Delta
}(X_{1},\ldots ,X_{n}),\]
as claimed.
\end{proof}
Let now $G$ be a discrete group, $S=\{g_{1},\ldots ,g_{n}\}$ a finite
symmetric set of generators (so that if $g\in S$, then $g^{-1}\in S$.)
Let $\lambda,\rho :G\to B(\ell ^{2}(G))$ be the left and right regular
representations given by $\lambda_g (f)(h)=f(g^{-1}h)$ and
$\rho_g(f)(h)=f(hg)$.  Then $J\rho_{g}J=\lambda_{g^{-1}}$.  Let $M$  be
the von Neumann algebra of $G$.

The following lemma is standard:
\begin{lem} \label{lem:groupvsM}
Consider the map $\phi: \ell^\infty(G)^{\oplus n} \to B(\ell^2(G))^n$
given by
$$\phi(\xi_1,\ldots,\xi_n)=(m_{\xi_1},\ldots,m_{\xi_n}),$$ where $m_f$
denotes
the operator of pointwise multiplication by $f\in\ell^\infty(G)$.  Then
\\ (a) $\phi(\ell^2(G)^{\oplus n}) \subset HS^n$; \\
\\
(b) For any closed $G$-invariant subspace $V\subset \ell^2(G)^{\oplus
n}$,
one has $$\dim_{M\bar\otimes M^o} \overline {M\phi(V)M} = \dim_{M} V.$$
\end{lem}

\begin{proof}
Part (a) is clear.

   For part (b), notice that we can identify $M=L(G)$
with $M^o$, and also $HS$ with
$\ell^2(G)\bar\otimes \ell^2(G)=\ell^2(G\times G)$.

With these identifications, if $\xi = \sum_g a_{g} \delta_g
\in \ell^2(G)$, with $\delta_g$ denoting the delta function at $g$,
then $\phi(\xi) = \sum a_g \delta_{g\times g} \in \ell^2(G\times
G )$.
Hence
$\phi$ is exactly the continuous extension to $L^2$ of the
induction map
$$1\otimes \cdot: L(G)\to (L(G\times G))\otimes_{L(G)} L(G)
= L(G)\bar \otimes L(G),$$ corresponding to the diagonal inclusion of
$G$ into $G\times G$ (see \cite[Theorem 3.3]{luck:foundations1}).
Now (b) follows because induction preserves dimension
\cite[Theorem 3.3]{luck:foundations1}.
\end{proof}

Recall that $\delta ^{(2)}(G)$ was defined by
$\delta ^{(2)}(G)=\beta _{1}(G)-\beta _{0}(G)+1$.

\begin{lem}
\label{lemma:delbtaBarvsdelta2} Let $S=\{g_1,\ldots,g_n\}$  be a
symmetric generating set for $G$.  Let
$U_j = \lambda_{g_j} $.  Then $\underline{\Delta }(U_{1},\ldots
,U_{n})\geq \delta ^{(2)}(G)$, as defined in Lemma~\ref{lem:delta2}.
\end{lem}
\begin{proof}
Let $\mathcal{G}$ be the Cayley graph of $G$ with respect to~$S$.
Identify $C_{(2)}^1(\mathcal{G})$ with $\ell^2(G)^{\oplus n}$, by
identifying the $j$-th copy of $\ell ^{2}(G)$ with edges labeled
$g_{j}$.

For $f\in \ell^\infty(G)$, denote by $\delta_{j}f$ the $j$-th component
of
$\delta f$
in this decomposition. Thus \( \delta _{j}f=\rho _{g_{j}}(f)-f\).
Let\[
a(f)=(\delta _{1}f,\ldots ,\delta _{n}f).\]

To prove the inequality
$\underline{\Delta}\geq \delta^{(2)}$, we need to provide a lower
estimate
on the dimension of the bimodule $\overline{H_0(U_1,\ldots,U_n)}$, and
so
we need some way of constructing elements in $H_0(U_1,\ldots,U_n)$.
In order to do that, we need some way of constructing bounded operators
$D$ so that $[D,J U_j J]\in HS$ for all $j=1,\ldots, n$.

We note that by Lemma \ref{lem:delta2}, we have that
$$\delta^{(2)}(G)=\dim_{L(G)}\overline{\{
c\in C^1_{(2)}(\mathcal{G}):c=\delta f\textrm{ for some $f\in
\ell^\infty(G)$}
\}}.$$
Now let $f\in \ell^\infty(G)$ be such that
$\delta f\in C^1_{(2)}(\mathcal{G})$. This is the same as saying that
$\delta_j(f)=
\rho_{g_j}(f)-f \in \ell^2(G)$ for each $j=1,\ldots,n$.

Denoting again by
$m_f$ the operator of multiplication by $f$, we have:
\begin{multline} \label{eq:groupvsM}
[m_{f},JU_{j}J]=m_{f}JU_{j}J-JU_{j}Jm_{f}=
JU_{j}J(JU_{j}^{-1}Jm_{f}JU_{j}J-m_{f})\\
= JU_jJ(m_{\rho_{g_j}(f)}-m_f)
=JU_{j}Jm_{\delta _{j}(f)}.
\end{multline}
Since $\delta _{j}(f)\in \ell ^{2}(G)$, we have that
$m_{\delta_j(f)}\in HS$ and so also $[m_{f},JU_{j}J]\in HS$.  Thus $m_f$
is a bounded operator whose commutators with
$J U_j J$, $j=1,\ldots,n$, are Hilbert-Schmidt operators.

Thus \[
A=\{([m_{f},JU_{1}J],\ldots ,[m_{f},JU_{n}J]):f\in \ell ^{\infty
}(G)\textrm{ s.t. }\delta f\in C_{(2)}^{1}(\mathcal{G})\}\subset
H_{0}(U_{1},\ldots ,U_{n}).\]
Since $H_{0}(U_{1},\ldots ,U_{n})$ is an $M,M$-bimodule,
   it will
suffice to prove that\[
\dim _{M\bar{\otimes }M^{o}}\overline{M\  A\  M}\geq \delta ^{(2)}(G)\]
(we'll actually prove that $\dim_{M\bar\otimes M^o} \overline{MAM} =
\delta^{(2)}(G)$.)

We now aim to use Lemma~\ref{lem:groupvsM} and the map $\phi$ defined
there.
Consider the $M,M$-bimodule isomorphism of $HS^n$ given by \[
\Psi :(\Xi _{1},\ldots ,\Xi _{n})\mapsto (JU_{1}^{-1}J\Xi _{1},\ldots
,JU_{n}^{-1}J\Xi _{n}).\]

    Then if $f\in \ell^\infty(G)$
with $\delta_j(f)\in \ell^2(G)$, $j=1,\ldots,n$, we have that
$$\Psi ( [m_f,JU_1J],\ldots,[m_f,JU_nJ]) = (m_{\delta_1 (f)},\ldots,
m_{\delta_n(f)})=\phi(\delta_1(f),\ldots,\delta_n(f)).$$
Hence
\begin{multline*}
\dim _{M\bar{\otimes }M^{o}}\overline{M\  A\  M}
=
\dim _{M\bar{\otimes }M^{o}}\overline{\Psi (M\  A\  M)}
\\ = \dim_{M\bar\otimes M^{o}} \overline{M\phi(\{c\in
C^1_{(2)}(\mathcal{G}):
c=\delta f\textrm{ for some }f\in \ell^\infty(G)\})M} \\
= \dim_{M} \overline{\{ c\in C_{(2)}^1 (\mathcal{G}) : c=\delta f
\textrm{ for some } f\in \ell^\infty(G)\}}
=\delta^{(2)}(M),
\end{multline*}
using Lemma~\ref{lem:groupvsM} and Lemma \ref{lem:delta2} in the last
two equalities.
\end{proof}

For any algebra $A$
generated by a self-adjoint set
of operators $X_{1},\ldots ,X_{n}$ on some Hilbert space $H$, and a
tracial state on $A$ given by $\tau(X)=\langle X\xi,\xi\rangle$, for
some
fixed $\xi\in H$,
let \[
\Delta (X_{1},\ldots ,X_{n})=n-\dim _{M\bar{\otimes
}M^{o}}\overline{\{(T_{1},\ldots ,T_{n})\in FR^{n}:\sum
_{j}[T_{j},JX_{j}J]=0\}},\]
where $M=W^{*}(X_{1},\ldots ,X_{n})$ is the von Neumann algebra
generated by $A$,
$FR$ stands for finite-rank
operators on $L^{2}(W^{*}(X_{1},\ldots ,X_{n}))$, and the closure
is taken in the Hilbert-Schmidt norm. This quantity was introduced in
\cite{connes-shlyakht:l2betti} and is related
to $L^{2}$-homology of $A$. The appearance of $FR$ comes from the fact
that this is the analogue of the space of compactly supported functions
on the group, in the same way that $HS$ is the analogue of the space
of square-summable functions.  One has:
$$\Delta (X_{1},\ldots ,X_{n})=
\beta _{1}(X_{1},\ldots ,X_{n})-\beta _{0}(X_{1},\ldots ,X_{n})+1$$
(we refer to \cite{connes-shlyakht:l2betti} for a definition of these
Betti numbers).
By \cite{connes-shlyakht:l2betti} one always has the inequality $$
\underline{\Delta}(X_1,\ldots,X_n)\leq \Delta(X_1,\ldots,X_n).$$ We
sketch
the proof for completeness.  Let $D\in B(L^2(M))$ be such that
$S_j = [JX_jJ,D]\in HS$, $j=1,\ldots,n$.  Then if $T_j\in FR$ satisfy $$
\sum_j [T_j,JX_jJ] = 0,$$ we have
\newcommand{\Tr}{\operatorname{Tr}}
$$
0=\Tr(\sum_j[T_j,JX_jJ]^* D)
   = \sum_j \Tr(T_j^* [JX_jJ,D]))
     = \sum_j \Tr(T_j^* S_j).
$$
Thus $(T_1,\ldots,T_n) \perp (S_1,\ldots,S_n)$ in $HS^n$.  Hence
$$H_0(X_1,\ldots,X_n)\perp \{(T_1,\ldots,T_n)\in FR^n:\sum_j
[T_j,JX_jJ]=0\}.$$
Since the Murray-von Neumann dimension of $HS^n$ over $M\bar\otimes
M^o$ is
$n$, it follows that $\underline{\Delta}\leq \Delta$.

\begin{cor}
\label{cor:DeltavsDeltabar}Let $Y_{1},\ldots ,Y_{n}$ be a self-adjoint
set of generators of $\mathbb{C}G$. Then\[
\Delta (Y_{1},\ldots ,Y_{n})=\underline{\Delta }(Y_{1},\ldots
,Y_{n})=\delta ^{(2)}(G),\]
where $\delta ^{(2)}(G)=\beta _{1}(G)-\beta _{0}(G)+1$.
\end{cor}
\begin{proof}
Since both $\Delta $ and $\underline{\Delta }$ don't depend on the
choice of generators of $\mathbb{C}G$, we may as well assume that
$Y_{1}=U_{1},\ldots ,Y_{n}=U_{n}$ correspond to a symmetric family
of generators of $G$. We then have by \cite[Theorem
3.3(c)]{connes-shlyakht:l2betti}
and Lemma \ref{lemma:delbtaBarvsdelta2} that \[
\delta ^{(2)}(G)\geq \Delta (U_{1},\ldots ,U_{n})\geq \underline{\Delta
}(U_{1},\ldots ,U_{n})\geq \delta ^{(2)}(G),\]
which forces all inequalities to be equalities.
\end{proof}

\section{Computation of free entropy dimension.}

Let $G$ be a finitely generated discrete group, and choose
$Y_{1},\ldots ,Y_{n}\in \mathbb{C}G$
to be self-adjoint elements in group algebra of $G$ that generate it as
a complex algebra.
One could
for example take $Y_{2j}=\Re \lambda_{g_{j}}
=\frac{1}{2}(\lambda_{g_j}+\lambda_{g_j}^{-1})$,
$Y_{2j-1}=\Im \lambda_{g_{j}} = \frac{1}{2i} (\lambda_{g_j}
- \lambda_{g_j}^{-1})$,
$j=1,\ldots ,n$ for some generators $g_{1},\ldots ,g_{n}$ of $G$.

\begin{thm} Let $G$ be a finitely generated group.
Let $Y_{1},\ldots ,Y_{n}$ be any self-adjoint generators of the group
algebra $\mathbb{C}G$, equipped with the canonical group trace $\tau$.
Then\[
\delta ^{*}(Y_{1},\ldots ,Y_{n})=\delta ^{\star }(Y_{1},\ldots
,Y_{n})=\beta _{1}(G)-\beta _{0}(G)+1.\]
In particular, $\delta ^{*}$ is an invariant of the algebra generated
by $Y_{1},\ldots ,Y_{n}$, taken with its trace.
\end{thm}

\begin{proof}
By \cite[Corollary 2.12]{shlyakht:qdim}\begin{equation}
\delta ^{*}(Y_{1},\ldots ,Y_{n})\geq \underline{\Delta }(Y_{1},\ldots
,Y_{n}).\label{eq:deltastarvsDeltabar}\end{equation}

By \cite[Theorem 4.4 and Corollary
4.6]{connes-shlyakht:l2betti},\begin{equation}
\Delta (Y_{1},\ldots ,Y_{n})\geq \delta ^{\star }(Y_{1},\ldots ,Y_{n})
\geq \delta ^{*}(Y_{1},\ldots
,Y_{n}).\label{eq:Deltavsdeltastar}\end{equation}

Combining (\ref{eq:deltastarvsDeltabar}), (\ref{eq:Deltavsdeltastar})
and Corollary \ref{cor:DeltavsDeltabar}, we find that\[
\delta ^{(2)}(G)=\Delta (Y_{1},\ldots ,Y_{n})\geq \delta ^{\star
}(Y_{1},\ldots ,Y_{n})\geq \delta ^{*}(Y_{1},\ldots ,Y_{n})\geq
\underline{\Delta }(Y_{1},\ldots ,Y_{n})=\delta ^{(2)}(G),\]
as claimed.
\end{proof}

\bibliographystyle{amsalpha}

\providecommand{\bysame}{\leavevmode\hbox to3em{\hrulefill}\thinspace}

\medskip{}
\noindent I. Mineyev, \textsc{250 Altgeld Hall, Department of
Mathematics,
University of Illinois at Urbana-Champaign, 1409 West Green St, Urbana,
IL 61801}

\noindent \texttt{mineyev@math.uiuc.edu,
http://www.math.uiuc.edu/\textasciitilde{}mineyev/math}

\medskip{}
\noindent D. Shlyakhtenko, \textsc{Department of Mathematics, UCLA,
Los Angeles, CA 90095}

\noindent \texttt{shlyakht@math.ucla.edu,
http://www.math.ucla.edu/\textasciitilde{}shlyakht/}
\end{document}